\documentclass[12pt]{article}

\usepackage{amsmath, amssymb}
\usepackage{amsfonts, mathrsfs}
\usepackage[cp1251]{inputenc}
\usepackage[active]{srcltx}
\usepackage[ukrainian, russian, english]{babel}

\usepackage{amsfonts, amssymb, mathrsfs}

\begin{document}

\begin{center}
	\textbf{\large \sc Approximation by Fourier sums in classes of differentiable functions with high exponents of smoothness}\\[1ex]
	\textbf{A.~S. Serdyuk and I.~V.  Sokolenko}\\[0.4ex]
	Institute of Mathematics of NAS of Ukraine, Kyiv, Ukraine\\[0.4ex]
	\textit{serdyuk@imath.kiev.ua, sokol@imath.kiev.ua} 
\end{center}

\begin{abstract}
We find asymptotic equalities for the exact upper bounds of approximations by Fourier  sums  of Weyl-Nagy classes  $W^r_{\beta,p}, 1\le p\le\infty,$ for  rapidly growing exponents of smoothness $r$ $(r/n\rightarrow\infty)$ in the  uniform metric. We obtain similar estimates for approximations of the classes $W^r_{\beta,1}$ in  metrics of the spaces $L_p, 1\le p\le\infty$.
\end{abstract}

	
\vskip 5mm

Let   $C$ be the space of continuous $2\pi$--periodic functions  $f$, in which the norm is defined by the equality
$$
{\|f\|_{C}=\max\limits_{t}|f(t)|},
$$
 $L_{p}$,
$1\le  p<\infty$, be the space of $2\pi$--periodic functions $f$ summable to the power $p$ on  $[-\pi,\pi)$, in which the norm is given by the formula
$$
\|f\|_{L_p}=\|f\|_{p}=\bigg(\int\limits_{-\pi}^{\pi}|f(t)|^pdt\bigg)^{1/p},
$$
 and $L_{\infty}$ be the space of measurable and essentially bounded   $2\pi$--periodic functions  $f$ with the norm
$$
\|f\|_{L_\infty}=\|f\|_{\infty}=\mathop {\rm ess \sup}\limits_{t} |f(t)|.
$$

Further, let $W^r_{\beta,p}, 1\le p\le \infty,$ be the sets of  all $2\pi$-periodic functions $f,$ representable as convolutions of the form
\begin{equation}\label{1}
f(x)=\frac{a_0}{2}+\frac{1}{\pi}\int\limits_{-\pi}^{\pi}
\varphi(x-t) B_{r,\beta}(t)dt, \ \ \ a_0\in\mathbb R,
\end{equation}
where  $B_{r,\beta}(\cdot)$ are Weyl-Nagy kernels  of the form
\begin{equation}\label{2}
B_{r,\beta}(t)=\sum\limits_{k=1}^\infty k^{-r}\cos\left(kt-\frac{\beta\pi}2\right),\quad r>0,\quad \beta\in\mathbb{R},
\end{equation}
 and functions $\varphi$ satisfy the condition
\begin{equation}\label{3'}
\varphi\in U_p^0=\left\{\varphi\in L_p: \|\varphi\|_p\le1,\   \int\limits_{-\pi}^{\pi}\varphi(t) dt=0\right\}.
\end{equation}

The classes $W^r_{\beta,p}$  are called  as  Weyl-Nagy classes (see, e.g., \cite{Sz.-Nagy1938, Stechkin1956, Stepanets1987, Stepanets2005}).

If  $r\in\mathbb N$ and $ \beta=r,\ $ then the functions of the form  (\ref{2})  are the well-known Bernoulli kernels and the classes $W^r_{\beta,p}$ coincide with the well-known classes $W^r_{p}$, which consist of $2\pi$-periodic functions   with absolutely continuous derivatives up to $(r-1)$-th order inclusive and such that $\|f^{(r)}\|_p\le1$ and $f^{(r)}(x)=\varphi(x)$ for almost everywhere $x\in\mathbb{R},\ $  where $\varphi$ is the function from (\ref{1}).

For arbitrary $\mathfrak N\subset X$, where $X=C$ or $L_p, \ 1\le p\le\infty,$ we consider the quantity
\begin{equation}\label{3}
\varepsilon_n(\mathfrak N)_X=\sup\limits_{f\in \mathfrak N}\|f(\cdot)-{S}_{n-1}(f; \cdot)\|_X,
\end{equation}
where  $S_{n-1}(f; x)$ is the partial Fourier sum of order $n-1$ of the
function $f$.

In the case of Weyl-Nagy classes $W^r_{\beta,\infty}$ and $X=C$
for the exact upper bounds (\ref{3})  the following asymptotic estimate holds
\begin{equation}\label{4}
\varepsilon_{n}(W^r_{\beta,\infty})_{C}=\frac4{\pi^2}\frac{\ln n}{n^r}+O\left(\frac1{n^r}\right),\quad r>0,\quad\beta\in\mathbb R.
\end{equation}

For $r\in\mathbb{N}$ and $\beta=r$  this estimate was obtained by A.N.~Kolmogorov \cite{Kolmogorov1985}, for arbitrary $r>0$  by V.T.~Pinkevich \cite{Pinkevich1940} and S.M.~Nikol'skii \cite{Nikol'skii1941}. In the general case the estimate (\ref{4}) follows from results, which were obtained in the works of A.V.Efimov \cite{Efimov1960} and  S.A.~Telyakovskii \cite{Telyakovskii1961}.

It should be also noticed, that a similar asymptotic equality holds for the classes $W^r_{\beta,1}$ in the metric of the space $L_1$, namely
\begin{equation}\label{5}
\varepsilon_{n}(W^r_{\beta,1})_{L_1}=\frac4{\pi^2}\frac{\ln n}{n^r}+O\left(\frac1{n^r}\right),\quad r>0,\quad\beta\in\mathbb R,
\end{equation}
(see \cite{Nikol'skii1946, Stechkin_Telyakovskii1967}).

In these works the parameters $r$ and $\beta$ of the Weyl-Nagy classes  were assumed to be fixed, and the question about the dependence of the remainder term in the estimates (\ref{4}) or (\ref{5}) on these parameters was not considered.

The character of the dependence on $r$ and $\beta$  of the remainder term in estimate (\ref{4}) was  investigated by I.G.~Sokolov \cite{Sokolov1955}, S.G.~Selivanova \cite{Selivanova1955}, G.I.~Natanson \cite{Natanson1961}, S.A.~Telyakovskii \cite{Telyakovskii1968, Telyakovskii1989} and S.B.~Stechkin \cite{Stechkin1980}.

In the work of S.B.~Stechkin \cite{Stechkin1980} the asymptotic behavior, as $n\rightarrow\infty$ and $r\rightarrow\infty$, of the quantities  $\varepsilon_{n}(W^r_{\beta,\infty})_{C}$ was completely investigated. Namely, he proved that for arbitrary $r\ge1$ and $\beta\in\mathbb{R}$ the following equality takes place
\begin{equation}\label{6}
\varepsilon_{n}(W^r_{\beta,\infty})_{C}=\frac1{n^r}\left(\frac8{\pi^2}\mathbf{K}(e^{-r/n})+O(1)\frac1r\right),
\end{equation}
 where
\begin{equation}\label{7}
\mathbf{K}(q)=\int\limits_{0}^{\pi/2}\frac{dt}{\sqrt{1-q^2\sin^2t}}
\end{equation}
is a complete elliptic integral of the first kind, and $O(1)$ is a quantity uniformly bounded with respect to $n, r$ and $\beta$.

Moreover, S.B.~Stechkin \cite[theorem 4]{Stechkin1980}  proved that for rapidly growing $r$ the remainder in estimate (\ref{6}) can be improved. Namely, for arbitrary $r\ge n+1$ and $\beta\in\mathbb{R}$ the following equality holds:
\begin{equation}\label{8}
\varepsilon_{n}(W^r_{\beta,\infty})_{C}=
\frac1{n^r}\left(\frac4{\pi} +O(1)\left(1+\frac1n\right)^{-r} \right),
\end{equation}
where $O(1)$ is a quantity uniformly bounded with respect to $n, r$ and $\beta$. If $ r/n\rightarrow\infty$, then the estimate $(\ref{8})$ becomes the asymptotic equality.

It also follows from \cite{Stechkin1980} that for the quantities  $\varepsilon_{n}(W^r_{\beta,1})_{L_1}$  the analogous estimates to (\ref{6}) and (\ref{8}) take place. Namely, for  $r\ge 1$ and $\beta\in\mathbb{R}$ uniformly with respect to the all analyzed parameters the following estimate is true
\begin{equation}\label{9}
\varepsilon_{n}(W^r_{\beta,1})_{L_1}=\frac1{n^r}\left(\frac8{\pi^2}\mathbf{K}(e^{-r/n})+O(1)\frac1r\right),
\end{equation}
where $\mathbf{K}(q)$ is defined by (\ref{7}), and for $r\ge n+1$ and $\beta\in\mathbb{R}$ uniformly with respect to all analyzed parameters the following estimate holds
\begin{equation}\label{10}
 \varepsilon_{n}(W^r_{\beta,1})_{L_1}=\frac1{n^r}\left(\frac4{\pi} +O(1)\left(1+\frac1n\right)^{-r} \right).
\end{equation}

Telyakovskii  \cite{Telyakovskii1989} showed that the remainder in formulas (\ref{8}) and (\ref{10}) can be replaced by a smaller one, namely, write $O(1)(1+\frac2n)^{-r}$ instead of $O(1)(1+\frac1n)^{-r}$.

In this paper for arbitrary values $1\le p\le\infty$ we establish generalized analogs of estimates (\ref{8}) and (\ref{10}) for  quantities $\varepsilon_{n}(W^r_{\beta,p})_{C}$ and $\varepsilon_{n}(W^r_{\beta,1})_{L_p}$, respectively.
Namely, as a consequence of the main result (Theorem 1) it  follows that for $r\ge n+1$, $\beta\in\mathbb{R}$ and  $1\le p\le\infty$ the following estimates hold:
\begin{equation}\label{11}
	\varepsilon_{n}(W^r_{\beta,p})_{C}=
	\frac1{n^r}\left(\frac{\|\cos t\|_{p'}}{\pi} +O(1)\left(1+\frac1n\right)^{-r}\right),
	\end{equation}
\begin{equation}\label{12}
	\varepsilon_{n}(W^r_{\beta,1})_{L_p}=
	\frac1{n^r}\left(\frac{\|\cos t\|_{p}}{\pi} +O(1)\left(1+\frac1n\right)^{-r} \right),
	\end{equation}
where $1/p+1/{p'}=1$ and  $O(1)$ are quantities uniformly bounded in all analyzed parameters. The estimates (\ref{11}) and (\ref{12}) are the asymptotic equalities, as  $r/n\rightarrow\infty$.

The main results of this paper are statements for the sets of functions, which are  more general compared to the Weyl-Nagy classes, namely, for the Stepanets classes $L^\psi_{\bar{\beta},p}$ and $C^\psi_{\bar{\beta},p}$ \cite{Stepanets1987, Stepanets2005}.

Denote by $L^\psi_{\bar{\beta},p},\ 1\le p\le\infty,$ the set of all $2\pi$-periodic functions $f$, representable for almost all $x\in\mathbb{R}$ as convolutions of the form
\begin{equation}\label{12'}
f(x)=\frac{a_0}2+\frac1\pi\int\limits_{-\pi}^{\pi}\Psi_{\bar{\beta}}(x-t)\varphi(t)dt, \ \ a_0\in\mathbb{R},
\end{equation}
of functions $\varphi$, which satisfy the conditions (\ref{3'}), with kernels $\Psi_{\bar{\beta}}\in L_1$, which Fourier series have the form
\begin{equation}\label{12''}
S[\Psi_{\bar{\beta}}](t)=\sum_{k=1}^{\infty}\psi(k)\cos\left(kt-\frac{\beta_k\pi}{2}\right),
\end{equation}
where $\psi=\{\psi(k)\}_{k=1}^\infty$ and $\bar{\beta}=\{\beta_k\}_{k=1}^\infty$ are fixed sequences of real numbers. In the case of $\bar{\beta}$ is a stationary  sequence, i.e.  $\beta_k\equiv\beta,\ \beta\in\mathbb{R},$ the kernels $\Psi_{\bar{\beta}}$ of the form (\ref{12''}) are denoted by $\Psi_{\beta}$, and the classes $L^\psi_{\bar{\beta},p}$ are denoted by  $L^\psi_{{\beta},p}$.

In this paper we consider the classes $L^\psi_{\bar{\beta},p}$, which are generated by  kernels of the form (\ref{12''}) with coefficients $\psi(k)>0$ such that
\begin{equation}\label{12'''}
\sum_{k=1}^{\infty}\psi(k)<\infty.
\end{equation}
In this case $S[\Psi_{\bar{\beta}}]=\Psi_{\bar{\beta}}\in C$ and, therefore, the imbedding $L^\psi_{\bar{\beta},p}\subset L_\infty$ is true. Since under condition (\ref{12'''}) the convolution of kernel $\Psi_{\bar{\beta}}$ with arbitrary function $\varphi\in U_p^0$ is a continuous function, then the set of all  $f$, representable in the form (\ref{12'}) for all $x\in\mathbb{R}$, is denoted by $C^\psi_{\bar{\beta},p}$.

For $\psi(k)=k^{-r}, r>1,$  classes $C^\psi_{\bar{\beta},p}$ are denoted by  $W^r_{\bar{\beta},p}$. If $\beta_k\equiv\beta,\ \beta\in\mathbb{R},$ then the classes $W^r_{\bar{\beta},p}$ are the Weyl-Nagy classes $W^r_{{\beta},p}$.

The main result of this paper is the following statement.

\textbf{Theorem 1.} \textit{
	Let $\ 1 \le p \le \infty,\ n\in\mathbb{N}$ and $\ \bar{\beta}=\{\beta_k\}_{k=1}^\infty$ be an arbitary
sequence of real numbers. Then for $r\ge n+1$ the following estimates
hold:
	\begin{equation}\label{_1t1}
	\varepsilon_{n}(W^r_{\bar\beta,p})_{C}=
	n^{-r}\Bigg(\frac{\|\cos t\|_{p'}}{\pi} +O(1)  \bigg(1+\frac1n\bigg)^{-r} \Bigg), \ \ \ \frac1p+\frac1{p'}=1,
	\end{equation}
	and
	\begin{equation}\label{_1t2}
	\varepsilon_{n}(W^r_{\bar\beta,1})_{L_p}=
	n^{-r}\Bigg(\frac{\|\cos t\|_{p}}{\pi} +O(1)  \bigg(1+\frac1n\bigg)^{-r} \Bigg),
	\end{equation}
where $O(1)$ are quantities uniformly bounded in all analyzed parameters.
}

\textit{Proof.} According to Theorem 4 from the work \cite{Serdyuk_2005_8} and  Theorem 4 from the work \cite{Serdyuk_2005_10} for arbitrary $\ 1 \le p \le \infty,\ \psi(k)>0,\  \bar{\beta}=\{\beta_k\}_{k=1}^\infty\ $ and $\ n\in\mathbb{N}$ taking into account the convergence of the series  $\sum\limits_{k=1}^{\infty}\psi(k)$ the following estimates hold:
\begin{equation}\label{t1d1}
\varepsilon_{n}(C^\psi_{\bar\beta,p})_{C}=
\frac{\|\cos t\|_{p'}}\pi\psi(n)+O(1)\sum_{k=n+1}^\infty \psi(k), \ \ \ \frac1p+\frac1{p'}=1,
\end{equation}
and
\begin{equation}\label{t1d2}
\varepsilon_{n}(L^\psi_{\bar\beta,1})_{L_p}=\varepsilon_{n}(C^\psi_{\bar\beta,1})_{L_p}=
\frac{\|\cos t\|_{p}}{\pi}\psi(n)+O(1)\sum_{k=n+1}^\infty \psi(k),
\end{equation}
where $O(1)$ are quantities uniformly bounded in all analyzed parameters.
We notice that   instead of the condition (\ref{12'''}) mentioned theorems have a stronger condition
\begin{equation}\label{t1d3}
\lim\limits_{k\rightarrow\infty} \frac{\psi(k+1)}{\psi(k)}=0.
\end{equation}
But, despite of it,  the proof of (\ref{t1d1}) and (\ref{t1d2}) does not require (\ref{t1d3}). The condition (\ref{t1d3}) provides the implementation of the relation
\begin{equation}\label{t1d4}
\sum\limits_{k=n+1}^\infty \psi(k)=o(1)\psi(n),
\end{equation}
and, therefore, guaranties that the estimates (\ref{t1d1}) and (\ref{t1d2})  are asymptotic equalities as $n\rightarrow\infty$.

Let put $\psi(k)=k^{-r}, r>1.$ Then, as noted above, $C^\psi_{\bar{\beta},p}=W^r_{\bar{\beta},p}, \ 1\le p\le\infty,$ and by virtue of (\ref{t1d1})
\begin{equation}\label{t1d5}
\varepsilon_{n}(W^r_{\bar{\beta},p})_{C}=
\frac{\|\cos t\|_{p'}}\pi\frac1{n^r}+O(1)\sum_{k=n+1}^\infty \frac1{k^r}, \ \ \ \frac1p+\frac1{p'}=1.
\end{equation}
Since  for arbitrary $n\in\mathbb{N},\ r\ge n+1$
\begin{equation}\label{t1d6}
\sum\limits_{k=n+1}^\infty\frac1{k^r}<\frac1{(n+1)^r}+\int\limits_{n+1}^{\infty}\frac{dt}{t^r}= \frac1{(n+1)^r}+ \frac1{(r-1)(n+1)^{r-1}}=\frac1{(n+1)^r}\frac{r+n}{r-1}\le
$$
$$
\le\frac1{n^r}\left(1+\frac1n\right)^{-r}\frac{2r-1}{r-1}\le\frac1{n^r}\left(1+\frac1n\right)^{-r}\left(2+\frac1n\right)\le 3 n^{-r}\left(1+\frac1n\right)^{-r},
\end{equation}	
then the estimate (\ref{_1t1}) follows from (\ref{t1d5}).

Similarly, in order to obtain the estimate  (\ref{_1t2}) we need to put $\psi(k)=k^{-r}, r>1,$ taking into account the equality  $C^\psi_{\bar{\beta},1}=W^r_{\bar{\beta},1}$ and estimate (\ref{t1d2}).  We get that
\begin{equation}\label{t1d7}
\varepsilon_{n}(W^r_{\bar\beta,1})_{L_p}=
\frac{\|\cos t\|_{p}}{\pi}\frac1{n^r}+O(1)\sum_{k=n+1}^\infty \frac1{k^r}.
\end{equation}

From formulas (\ref{t1d6}) and (\ref{t1d7}) we obtain (\ref{_1t2}). Theorem 1 is proved.

\textbf{Remark 1.} \textit{If the conditions of Theorem 1 are satisfied and, in addition,
\begin{equation}\label{1z1}
\frac{r}{n}\rightarrow\infty,
\end{equation}
then formulas (\ref{_1t1}) and (\ref{_1t2}) are asymptotic equalities.}

Indeed, since the sequence $\left(1+1/n\right)^{n+1}$ is monotonically decreasing to number $e$, we have that
\begin{equation}\label{1z2}
\left(1+\frac1n\right)^{-r}=\left(\left(1+\frac1n\right)^{n+1}\right)^{-\frac r{n+1}}
\le e^{-\frac r{n+1}}.
\end{equation}
As it follows from (\ref{1z2}) under condition (\ref{1z1}) $
 \left(1+\frac{1}{n}\right)^{-r}\rightarrow 0,
$
and, therefore, the estimates (\ref{_1t1}) and (\ref{_1t2}) are asymptotic equalities.

In the case of $\bar{\beta}=\{\beta_k\}_{k=1}^\infty$ are the stationary consequences, i.e. $\beta_k\equiv\beta,\ \beta\in\mathbb{R},$ the formulas  (\ref{11}) and (\ref{12}) follow from (\ref{_1t1}) and (\ref{_1t2}), respectively.

The following statement follows from Theorem 1 in the case $p=2$.

\textbf{Corollary 1.} \textit{
	Let $\ \bar{\beta}=\{\beta_k\}_{k=1}^\infty$ be an arbitrary
sequence of real numbers, $\ n\in\mathbb{N}$ and $r\ge n+1$. Then
	\begin{equation}\label{1c1}
\varepsilon_{n}(W^r_{\bar\beta,2})_{C}=\varepsilon_{n}(W^r_{\bar\beta,1})_{L_2}=
	\frac1{n^r}\left(\frac1{\sqrt\pi} +O(1)  \bigg(1+\frac1n\bigg)^{-r} \right),
	\end{equation}
	where $O(1)$ is quantity uniformly bounded in all analyzed parameters.
}

In order to prove the first equality in (\ref{1c1}) it is enough to use the results of the works \cite{Serdyuk_Sokolenko2011} and \cite{Serdyuk_Sokolenko2013}. As it follows from these results for all $r>\frac12$ and $n\in\mathbb{N}$
\begin{equation}\label{1c2}
\varepsilon_{n}(W^r_{\bar\beta,2})_{C}=\varepsilon_{n}(W^r_{\bar\beta,1})_{L_2}
=\frac1{\sqrt\pi}\left(\sum\limits_{k=n}^\infty\frac1{k^{2r}}\right)^{1/2}.
\end{equation}

To prove the second equality in (\ref{1c1}) it is sufficient to use the estimates (\ref{_1t1}) and (\ref{_1t2}) setting $p=p'=2$.
However, the estimate (\ref{1c1}) can also be obtained directly from equation (\ref{1c2}),
taking into account that
\begin{equation}\label{1c3}
\frac1{\sqrt\pi}\left(\sum\limits_{k=n}^\infty\frac1{k^{2r}}\right)^{1/2}=\frac1{\sqrt\pi} \left(\frac1{n^r}+O(1)\sum\limits_{k=n+1}^\infty\frac1{k^r}\right),
\end{equation}
and using relation (\ref{t1d6}) for $r\ge n+1$.

We notice, since the Hurwitz zeta function $\zeta(s,l)=\sum\limits_{m=0}^\infty (l+m)^{-s}, Re(s) > 1$,  $Re(l) > 0,$ has an integral representation
$$
\zeta(s,l)=\frac1{\Gamma(s)}\int\limits_0^\infty\frac{t^{s-1}e^{-lt}}{1-e^{-t}}dt,
$$
then the formula (\ref{1c2}) can be rewritten in the equivalent form
$$
\varepsilon_{n}(W^r_{\bar\beta,2})_{C}=\varepsilon_{n}(W^r_{\bar\beta,1})_{L_2}=
\frac1{\sqrt{\pi\Gamma(2r)}}\left(\int\limits_0^\infty\frac{t^{2r-1}e^{-nt}}{1-e^{-t}}dt\right)^{\frac12}.
$$

\end{document}